\documentclass[10pt]{amsart}
\usepackage{graphicx}
\usepackage{latexsym}
\usepackage{fancyhdr}
\usepackage{amsmath, amssymb}
 \usepackage[utf8]{inputenc}
\usepackage[all]{xy}
\usepackage{float}
\usepackage{tikz}
\usepackage{pdflscape}
\usepackage{longtable}
\usepackage{rotating}
\usepackage{hyperref}
\usepackage{subfigure}
\usepackage{tensor}
\usepackage{caption}
\usepackage{enumerate}

\usepackage{color}
\theoremstyle{plain}
\newtheorem{thm}{Theorem}

\newtheorem{lem}[thm]{Lemma}

\theoremstyle{definition}
\newtheorem{defi}[thm]{Definition}
\theoremstyle{remark}
\newtheorem{rem}[thm]{Remark}
\numberwithin{equation}{section}

\newcommand{\lgw}{\longrightarrow}

\newcommand{\ovl}{\overline}

\newcommand{\G}{\Gamma}

\renewcommand{\S}{\mathcal S}

\newcommand{\Z}{\mathbb{Z}}

\newcommand{\Graph}{\operatorname{Graph}}

\newcommand{\N}{\mathbb{N}}

\newcommand{\C}{\mathbb{C}}

\renewcommand{\phi}{\varphi}

\newcommand{\e}{\varepsilon}

\let\mathscr\mathcal

\DeclareMathOperator{\re}{Re}
\DeclareMathOperator{\im}{Im}

\begin{document}
\title[Transcendental holomorphic maps between  algebraic manifolds]{Transcendental  holomorphic maps between real algebraic manifolds in a complex space}
\author{Guillaume Rond}
\email{guillaume.rond@univ-amu.fr}
\address{Aix-Marseille Universit\'e, LASOL, UMI2001, UNAM, Mexico}

\begin{abstract}
We give an example of a real algebraic manifold embedded in a complex space that does not satisfy the Nash-Artin approximation Property. This Nash-Artin approximation Property is closely related to the problem of determining when the biholomorphic equivalence for germs of real algebraic manifolds coincides with the algebraic equivalence. This example is an elliptic Bishop surface, and its construction is based on the functional equation satisfied by the generating series of some walks restricted to the quarter plane.
\end{abstract}

\thanks{The  author is deeply grateful to the UMI LASOL of the CNRS where this project has been carried out.}


\maketitle

\section{Introduction}
An important problem in complex geometry is to classify germs of real analytic manifolds in  $\C^N$ up to biholomorphic equivalence (see \cite{BER} for an introduction to this problem). This problem goes back to Poincar\'e \cite{Po}, and E. Cartan, for germs of real analytic smooth hypersurfaces in $\C^2$ \cite{Ca}, then S. S. Chern and J. K. Moser, for germs of real analytic smooth hypersurfaces in $\C^N$ for $N\geq 2$ \cite{CM}, gave a complete description of this classification, when the hypersurfaces are assumed to be Levi nondegenerate. More precisely, S. S. Chern and J. K. Moser first gave a complete classification up to formal biholomorphisms, 
and then they showed that any formal biholomorphism between Levi nondegenerate real hypersurfaces is necessarily convergent.\\
Therefore a natural question was to understand when the biholomorphic equivalence and the formal (biholomorphic) equivalence coincide. This question has been widely studied since the work \cite{CM}; the reader can consult \cite{Mir} for a general presentation of this problem. The first negative answer to this question has been given in \cite{MW}: the authors considered a particular example of a germ of a real algebraic smooth surface $(M,0)$ in $\C^2$ for which there exists a germ of a real analytic smooth surface $(M',0)\subset (\C^2,0)$, such that $(M,0)$ and $(M',0)$ are formally equivalent but not biholomorphically equivalent. This surface $M$ has the following particular property: its tangent space at  any point near the origin is totally real, but its tangent space at the origin is a complex line (in other words it has a CR-singularity at the origin). Such a surface is called a Bishop surface. \\
On the other hand, this question is known to have a positive answer for several examples of CR-manifolds (see \cite{Mir} for an introduction to the subject), and for quite a while, it remained open whether, for any germ of a real analytic CR-manifold $(M,0)$, the formal equivalence class of $(M,0)$ was the same as its biholomorphic equivalence class. This question has been recently answered in the negative in \cite{KS}.\\
In the case of real algebraic (not necessarily CR) manifolds in $\C^N$, one can define the notion of algebraic (biholomorphic) equivalence. One says that two germs of real algebraic manifolds $(M,0)$, $(M',0)$ in $\C^N$ are algebraically equivalent if there is a germ of biholomorphic map $h:(\C^N,0)\lgw (\C^N,0)$ such that $h(M)=M'$, and such that the components of $h(z)$ are algebraic, that is, they are given by algebraic power series. The question of whether biholomorphic equivalence implies  algebraic equivalence of germs of real algebraic manifolds has first been asked in \cite{BER}. This is even stated as a conjecture in \cite{BMR}. Up to now the answer to this question remains completely open (see also \cite{Mir}), even if it is known to be true in some important cases (as the case of hypersurfaces \cite{BMR} - see also \cite{Z} and \cite{KLS}).\\
A related problem has been introduced by N. Mir in \cite{Mi2}. He introduced the  notion of Nash-Artin approximation Property for a smooth real algebraic set $M$  (see Definition \ref{Nash-Artin} below). The Nash-Artin approximation Property for $M$ implies that the biholomorphic equivalence class of $(M,0)$ coincides with the algebraic equivalence class of $(M,0)$. In fact this approximation property is natural since, for a lot of the cases where these two notions of equivalence are known to coincide, the Nash-Artin property is also satisfied. Therefore, a natural question is whether  this approximation property is always satisfied for real algebraic manifolds. In the case of real algebraic CR-manifolds, this is even a conjecture (see \cite{Mir}).\\
Up to now, there were no known examples of real algebraic manifolds not satisfying this Nash-Artin approximation Property.
The aim of this note is to provide such examples. These examples are constructed from the functional equations satisfied by the generating series of some walks restricted to the quarter plane. These examples are in fact elliptic Bishop surfaces.\\
\\
The author would like to thank Nordine Mir for his comments and for pointing out that the surfaces $M_\S$ are elliptic Bishop surfaces.

\section{The Nash-Artin approximation Property}\label{NA}

\begin{defi}\label{Nash-Artin}\cite{Mi2}
Let $(M,0)\subset (\C^N,0)$ be the germ of a real algebraic manifold. We say that $(M,0)$ has the Nash-Artin approximation Property if
\begin{itemize}
\item[i)] for every $N'$ and every real algebraic set $\G\subset \C^N\times\C^{N'}$
\item[ii)] for every germ of holomorphic map $h:(\C^N,0)\lgw \C^{N'}$
\item[iii)] for every $c\in\N$
\end{itemize}
such that
$$\Graph(h)\cap(M\times\C^{N'})\subset \G,$$
there is a holomorphic algebraic map germ $h_c:(\C^N,0)\lgw \C^{N'}$ such that

$$\Graph(h_c)\cap(M\times\C^{N'})\subset \G,$$
and $h_c$ agrees with $h$ up to order $c$.

\end{defi}
The terminology comes from the Artin approximation Theorem \cite{Ar}.
\begin{rem}
Assume that the ideal defining $M$ is generated by the real valued polynomials $f_1(z,\ovl z)$, \ldots, $f_p(z,\ovl z)$ and the ideal defining $\G$ is generated by the real valued polynomials $g_1(z,w,\ovl z,\ovl w)$, \ldots, $g_s(z,w,\ovl z,\ovl w)$. Here $z=(z_1,\ldots, z_N)$ and $w=(w_1,\ldots, w_{N'})$.
Then the condition 
$$\Graph(h)\cap(M\times\C^{N'})\subset \G,$$
is equivalent to the existence of series $k_{j,l}(z,\ovl z)$, for $1\leq j\leq s$ and $1\leq l\leq p$,  such that

\begin{equation}\label{key_nash} 
g_j(z,h(z),\ovl z,\ovl{h(z)})+\sum_{l=1}^pk_{j,l}(z,\ovl z)f_l(z,\ovl z)=0\ \ \forall j=1,\ldots, s.
\end{equation}

\end{rem}

\section{Generating series of walks restricted to the quarter plane}
We consider the following situation: We fix a set of \emph{steps} $\mathcal S\subset \{-1,0,1\}^2\setminus\{(0,0)\}$. For all $(i,j)\in\N^2$ and $n\in\N$, we denote by $a_{i,j,n}$ the number of walks with steps in $\mathcal S$ of length $n$, starting at the origin and ending at $(i,j)$, and remaining in the quarter plane $\N^2$. The associated generating series is defined by
$$Q(x,y,t):=\sum_{i,j,n\in\N}a_{i,j,n}x^iy^jt^n\in \Z[[x,y,t]].$$
The reader may consult \cite{BM} for a general account of the study of these generating series. 
We recall here  that $Q(x,y,t)$ is the solution of the equation
\begin{equation}\begin{split}\label{func_eq} xy=\left(xy-t\sum_{(a,b)\in\mathcal S}x^{a+1}y^{b+1}\right)Q(x,y,t)&+ty\sum_{j\mid (-1,j)\in \mathcal S}y^jQ(0,y,t)+\\
+tx&\sum_{i\mid (i,-1)\in \mathcal S}x^iQ(x,0,t)+\e tQ(0,0,t)\end{split}\end{equation}
where $\e=1$ if $(-1,-1)\in \mathcal S$ and $\e=0$ otherwise.\\
The term $xy-t\sum_{(a,b)\in\mathcal S}x^{a+1}y^{b+1}$ is called the kernel of the equation and is denoted by $K_\S(x,y,t)$. In fact, by the division theorem of formal power series (see for instance \cite[Example 1.14]{Ro}), the equation
\begin{equation}\label{func2}xy=K_\S(x,y,t)l(x,y,t)+r(x,y,t)\end{equation}
has a unique solution  $(l(x,y,t),r(x,y,t))\in\C[[x,y,t]]$ such that none of the monomials of $r(x,y,t)$ are divisble by $xy$. In particular, $r(x,y,t)$ can be written as
$r(x,y,t)=h(x,t)+g(y,t)$ where $(h(x,t), g(y,t))\in\C[[x,t]]\times\C[[y,t]]$. Therefore we have necessarily 
$l(x,y,t)=Q(x,y,t)$ and
$$ r(x,y,t)=ty\sum_{j\mid (-1,j)\in \mathcal S}y^jQ(0,y,t)
+tx\sum_{i\mid (i,-1)\in \mathcal S}x^iQ(x,0,t)+\e tQ(0,0,t).$$
We have the following lemma (for a general account of Bishop surfaces, see  \cite{Bi}, \cite{Mo} or \cite{HY} for example) :

\begin{lem}\label{lem_bishop}
We denote by $M_\S$ the real algebraic surface defined by $K_\S(w,\ovl w,z)=0$. Then $M_S$ is smooth at the origin if and only if $(-1,-1)\in\S$. In this case the germ $(M_S,0)$ is the germ of a Bishop surface with a Bishop invariant equal to 0.
\end{lem}

\begin{proof}
We have 
$$\frac{\partial K_\S}{\partial w}(w,\ovl w,z)=\ovl w-z\sum_{(a,b)\in\mathcal S}(a+1)w^{a}\ovl w^{b+1},$$
$$\frac{\partial K_\S}{\partial \ovl w}(w,\ovl w,z)=w-z\sum_{(a,b)\in\mathcal S}(b+1)w^{a+1}\ovl w^{b},$$
$$\frac{\partial K_\S}{\partial z}(w,\ovl w,z)=\sum_{(a,b)\in\mathcal S}w^{a+1}\ovl w^{b+1}.$$
Thus, $\frac{\partial K_\S}{\partial w}(0,0,0)=\frac{\partial K_\S}{\partial \ovl w}(0,0,0)=0$. Therefore $M_\S$ is smooth at the origin if and only if 
$\frac{\partial K_\S}{\partial z}(0,0,0)\neq 0$, that is, if and only if $(-1,-1)\in\S$.\\
In this case, the germ $(M_\S,0)$ is also defined by
$$\phi(w,\ovl w, z):=z-w\ovl w\left(1+\sum_{(a,b)\in\mathcal S\setminus\{(-1,-1)\}}w^{a+1}\ovl w^{b+1}\right)^{-1}=0.$$
But $\phi$ can be expanded as
$$\phi(w,\ovl w,z)=z-w\ovl w+E(w,\ovl w)$$
where $E(w,\ovl w)$ is the germ of a real analytic function  vanishing to order  at least three at the origin. Therefore $(M_S,0)$ is the germ of a Bishop surface with a Bishop invariant equal to 0.
\end{proof}

To such a kernel $K_\S$ is usually associated a group  of birational automorphisms $G(\S)$ that preserves $K_\S$. By \cite[Theorem 1]{KR}, this group is finite if and only if the series $Q_\S$ is D-finite (or holonomic). In fact they show the following stronger result: if $G(\S)$ is infinite then $Q(x,0,t)$ and $Q(0,y,t)$ are not D-finite. In particular, in this case,  $Q(x,0,t)$ and $Q(0,y,t)$ are transcendental convergent power series. There are 56 such walks and their list  can be found in \cite[Table 4]{BM} for instance. \\
Among these 56 walks with infinite group, only 7 satisfy Lemma \ref{lem_bishop} while having a Kernel satisfying the following symmetry:
\begin{equation}\label{sym} K_\S(w,\ovl w,z)=K_\S(\ovl w, w, z).\end{equation}
These correspond to the following sets $\S$:
\begin{figure}[H]\label{fig}

\fbox{\begin{tikzpicture}[scale=0.6]

\draw[line width=0.5pt] (0,0) -- (1,0);
\draw[line width=0.5pt] (0,0) -- (0,1);
\draw[line width=0.5pt] (0,0) -- (1,1);
\draw[line width=0.5pt] (0,0) -- (-1,-1);

\end{tikzpicture}\hspace{0.5cm}
\begin{tikzpicture}[scale=0.6]

\draw[line width=0.5pt] (0,0) -- (-1,0);
\draw[line width=0.5pt] (0,0) -- (0,-1);
\draw[line width=0.5pt] (0,0) -- (1,1);
\draw[line width=0.5pt] (0,0) -- (-1,-1);

\end{tikzpicture}\hspace{0.5cm}
\begin{tikzpicture}[scale=0.6]

\draw[line width=0.5pt] (0,0) -- (-1,0);
\draw[line width=0.5pt] (0,0) -- (0,-1);
\draw[line width=0.5pt] (0,0) -- (0,1);
\draw[line width=0.5pt] (0,0) -- (1,0);
\draw[line width=0.5pt] (0,0) -- (-1,-1);

\end{tikzpicture}\hspace{0.5cm}
\begin{tikzpicture}[scale=0.6]

\draw[line width=0.5pt] (0,0) -- (-1,1);
\draw[line width=0.5pt] (0,0) -- (1,-1);
\draw[line width=0.5pt] (0,0) -- (0,1);
\draw[line width=0.5pt] (0,0) -- (1,0);
\draw[line width=0.5pt] (0,0) -- (-1,-1);

\end{tikzpicture}\hspace{0.5cm}
\begin{tikzpicture}[scale=0.6]

\draw[line width=0.5pt] (0,0) -- (-1,1);
\draw[line width=0.5pt] (0,0) -- (1,-1);
\draw[line width=0.5pt] (0,0) -- (0,1);
\draw[line width=0.5pt] (0,0) -- (1,0);
\draw[line width=0.5pt] (0,0) -- (-1,-1);
\draw[line width=0.5pt] (0,0) -- (1,1);

\end{tikzpicture}\hspace{0.5cm}
\begin{tikzpicture}[scale=0.6]

\draw[line width=0.5pt] (0,0) -- (-1,1);
\draw[line width=0.5pt] (0,0) -- (1,-1);
\draw[line width=0.5pt] (0,0) -- (0,-1);
\draw[line width=0.5pt] (0,0) -- (-1,0);
\draw[line width=0.5pt] (0,0) -- (-1,-1);
\draw[line width=0.5pt] (0,0) -- (1,1);

\end{tikzpicture}\hspace{0.5cm}
\begin{tikzpicture}[scale=0.6]

\draw[line width=0.5pt] (0,0) -- (-1,1);
\draw[line width=0.5pt] (0,0) -- (1,-1);
\draw[line width=0.5pt] (0,0) -- (0,-1);
\draw[line width=0.5pt] (0,0) -- (-1,0);
\draw[line width=0.5pt] (0,0) -- (-1,-1);
\draw[line width=0.5pt] (0,0) -- (0,1);
\draw[line width=0.5pt] (0,0) -- (1,0);
\end{tikzpicture}}

\caption{}
\label{fig}
\end{figure}
For these walks, by symmetry, the solutions $(l,r)$ of \eqref{func2} satisfy the relation $r(x,y,t)=r(y,x,t)$ and have positive integer coefficients.\\
In fact, it has been recently shown that the generating series corresponding to these 7 types of walks are among those that  are $x$-hypertranscendental and $y$-hypertranscendental \cite{DHRS}. That is, $Q(x,0,t)$ (resp. $Q(0,y,t)$) is not a solution, as a function of $x$ (resp. of $y$), of a polynomial differential equation. Therefore, $Q(x,y,t)$ is not a solution of a polynomial differential equation, and neither is $r(x,y,t)$.

\section{Example}
Let $M_\S\subset \C^2_{z,w}$ be the real algebraic manifold defined by $K_\S$ where $\S$ is one of the sets given in Figure \ref{fig}. 
Since
$$\sum_{(a,b)\in\mathcal S}w^{a+1}\ovl w^{b+1}=\sum_{(a,b)\in\mathcal S}\ovl w^{a+1}w^{b+1}$$
by \eqref{sym}, we easily check that $M$ is defined by the real algebraic equations:
$$\left\{\begin{array}{c}\im(z)=0\\
|w|^2-z\sum_{(a,b)\in\mathcal S}w^{a+1}\ovl w^{b+1}=0\end{array}\right..$$
Let $\G$ be the subset of $\C^2_{z,w}\times\C_{h}$ defined by
$$|w|^2+h+\ovl h=0$$
or, equivalently, by
$$ |w|^2+2\re(h) = 0.$$
Then $\G$ is a real algebraic smooth hypersurface.
Now let 
$$h:(\C^2,0)\lgw (\C,0)$$
be a germ of a (formal) holomorphic map such that
$\Graph(h)\cap(M_\S\times \C)\subset\G.$
By \eqref{key_nash}, this is equivalent to saying that there exist power series $k(z,w,\ovl z,\ovl w)$ and $l(z,w,\ovl z,\ovl w)$ such that the following relations are satisfied:

$$|w|^2+2\re(h(z,w))
+yk(z,w,\ovl z,\ovl w)+K_\S(w,\ovl w, z)l(z,w,\ovl z,\ovl w)=0.$$
 Equivalently we have
\begin{equation}\label{main_eq1}w\ovl w+h(x,w)+\ovl h(x,\ovl w)+K_\S(w,\ovl w, x)l(x,w,\ovl w)=0.\end{equation}
for some $l(x,w,\ovl w)\in\C[[x,w,\ovl w]]$.
But \eqref{main_eq1} is exactly \eqref{func2} whose unique solution is  the generating series $l$ that counts the number of walks restricted to the quarter plane by the length and by the end point, and whose set of elementary steps is $\mathcal S$. In this case we have
$$h(x,w)+\ovl h(x,\ovl w)=r(w,\ovl w,x)$$ is a 
 transcendental power series (with real coefficients) as explained before.  \\
This proves that the  germs of holomorphic maps $h$ such that 
$$\Graph(h)\cap (M_\S\times\C)\subset \G$$ are  not algebraic.\\
This shows that $M_\S$ does not satisfy the Nash-Artin approximation Property.




\begin{thebibliography}{00}

\bibitem[Ar69]{Ar} M. Artin, Algebraic approximation of structures over complete local rings, \textit{Publ. Math. IHES}, \textbf{36}, (1969), 23-58.
 

\bibitem[BER00]{BER} M. S. Baouendi, P. Ebenfeld, L. P. Rothschild, Local geometric properties of real submanifolds in complex space, \emph{Bull. AMS}, \textbf{37}, no. 3, (2000), 309-336.


\bibitem[BMR02]{BMR}  M. S. Baouendi, N.  Mir, L. P.  Rothschild,  Reflection ideals and mappings between generic submanifolds in complex space, \emph{J. Geom. Anal.}, \textbf{12}, No. 4, (2002), 543-580. 




\bibitem[Bi65]{Bi} E. Bishop, Differentiable manifolds in complex Euclidean space,
\emph{Duke Math. J.}, \textbf{32},  (1965), 1-21. 

\bibitem[BMM09]{BM} M. Bousquet-M\'elou, M.  Mishna,  Walks with small steps in the quarter plane, \emph{Contemp. Math.}, \textbf{520},  (2010), 1-40.

\bibitem[Ca32]{Ca} \'E. Cartan, Sur la g\'eom\'etrie pseudo-conforme  des hypersurfaces de l'espace de deux variables complexes, \emph{Ann. Mat. Pura Appl.}, \textbf{11(4)}, (1932), 17-90.

\bibitem[CM75]{CM} S.S. Chern, J.K.  Moser,  Real hypersurfaces in complex manifolds, \emph{Acta Math.}, \textbf{133}, (1974), 219-271.

\bibitem[DHRS18]{DHRS} T. Dreyfus, C. Hardouin, J. Roques, M. Singer, 
On the nature of the generating series of walks in the quarter plane,
\emph{Invent. Math.}, \textbf{213}, No. 1,  (2018), 139-203. 

\bibitem[HY09]{HY} X. Huang, Xiaojun; W. Yin, 
A Bishop surface with a vanishing Bishop invariant, \emph{Invent. Math.}, \textbf{176}, No. 3, (2009), 461-520. 


\bibitem[KR12]{KR} I. Kurkova, K.  Raschel, On the functions counting walks with small steps in the quarter plane, \textit{Publ. Math. Inst. Hautes \'Etudes Sci.}, \textbf{116}, (2012), 69-114. 


\bibitem[KS16]{KS} I. Kossovskiy, R. Shafikov, Divergent CR-equivalences and meromorphic differential equations,
\emph{J. Eur. Math. Soc.}, \textbf{18}, (2016), no. 12, 2785-2819. 

\bibitem[KLS16]{KLS} I. Kossovskiy, B. Lamel, L. Stolovitch, Equivalence of Cauchy-Riemann manifolds and multisummability theory, ArXiv:1612.05020.



\bibitem[Mir12]{Mi2} N. Mir,  Algebraic approximation in CR geometry, \emph{J. Math. Pures Appl.}, \textbf{98}, (2012), 72-88.

\bibitem[Mir13]{Mir} N. Mir,  Artin  approximation theorems and Cauchy-Riemann geometry, \emph{Methods Appl. Anal.}, \textbf{21} (4), (2104), 481-502.

\bibitem[Mo85]{Mo} J. Moser,  Analytic surfaces in $\C^2$ and their local hull of holomorphy, 
\emph{Ann. Acad. Sci. Fenn., Ser. A I, Math.}, \textbf{10}, (1985), 397-410. 

\bibitem[MW83]{MW} J. K. Moser, S.  Webster, Normal forms for real surfaces in $\C^2$ near complex tangents and hyperbolic surface transformations, \emph{Acta Math.}, \textbf{150}, (1983), 255-296.

\bibitem[Po07]{Po} H. Poincar\'e, Les fonctions analytiques de deux variables et la repr\'esentation conforme, \emph{Rend. Circ. Mat. Palermo}, \textbf{II. Ser. 23}, (1907), 185-220. 


\bibitem[Ro18]{Ro} G. Rond, Artin approximation, \emph{J. Singul.}, \textbf{17}, (2018), 108-192.

\bibitem[Za99]{Z} D. Zaitsev, Algebraicity of local holomorphisms between real-algebraic submanifolds of complex spaces,
\emph{Acta Math.}, \textbf{183}, No. 2, (1999), 273-305. 



\end{thebibliography}
\end{document}